\mag=1000 \hsize=6.5 true in \vsize=8.7 true in
    \baselineskip=15pt
\vglue 1.5 cm

%%%%%%%%%%%%%%%%%%%%%%%%%%%%%%%%%%%%%%%%%%%%%%%%%%%%%%%%%%%%%%%%
    %This is an up-down, i.e. vertical contains

%%%%%%%%%%%%%%%%%%%%%%%%%%%%%%%%%%%%%%%%%%%%%%%%%%%%%%%%%%%%%%%%%%%%%%%%%
          % This macro makes an empty box to be used as "qed".

    \def\qed{$\rlap{$\sqcap$}\sqcup$}

%%%%%%%%%%%%%%%%%%%%%%%%%%%%%%%%%%%%%%%%%%%%%%%%%%%%%%%%%%%%%%%%%%%%%%%%%%
    % The following is what allows me to make gothic letters

\font\tengothic=eufm10 \font\sevengothic=eufm7
\newfam\gothicfam
        \textfont\gothicfam=\tengothic
        \scriptfont\gothicfam=\sevengothic
\def\goth#1{{\fam\gothicfam #1}}

%%%%%%%%%%%%%%%%%%%%%%%%%%%%%%%%%%%%%%%%%%%%%%%%%%%%%%%%%%%%%%%%%%%%%%%%%%%
    % The following allows the use of the family Bbb of double stroked
    % letters

     \font\tenmsb=msbm10              \font\sevenmsb=msbm7
%   \font\tenmsb=msym10              \font\sevenmsb=msym7
\newfam\msbfam
        \textfont\msbfam=\tenmsb
        \scriptfont\msbfam=\sevenmsb
\def\Bbb#1{{\fam\msbfam #1}}

%%%%%%  %%%%%%%%%%%%%%%%%%%%%%%%%%%%%%%%%%%%%%

\def\PP#1{{\Bbb P}^{#1}}
%\def\PP#1{{\bf P}^{#1}}
    %%%%%%%%%%%%%%%%%%%%%%%%%%%%%%
\def\PPP {{\Bbb P}}

%\def\PPP{{\bf P}}
%%%%%%%%%%%%%%%%%%%%%%%%%%%%%%%%%%
%for references

\def\ref#1{[{\bf #1}]}
    %%%%%%%%%%%%%%%%%%%%%%%%%%%%%

\def\PP#1{{\Bbb P}^{#1}}
%\def\PP#1{{\bf P}^{#1}}
   %%%%%%%%%%%%%%%%%%%%%%%%%%%%%%
\def\PPP {{\Bbb P}}
%\def\PPP{{\bf P}}
%%%%%%%%%%%%%%%%%%%%%%%%%%%%%%%%%%
%for references
\def\ref#1{[{\bf #1}]}
   %%%%%%%%%%%%%%%%%%%%%%%%%%%%%
%for primes

\def\a{\bigskip \par \noindent}
\def\b{\par \noindent}

   %%%%%%%%%%%%%%%%%%%%%%%%%%%%%%%%%%%%%%%%%%%%%%%%%%%%%%
\def\proof{\bigskip \par \noindent {\bf Proof.  }}

   %%%%%%%%%%%%%%%%%%%%%%%%%%%%%
\def\prop#1{\bigskip \par \noindent {\bf #1. Proposition.}}
   %%%%%%%%%%%%%%%%%%%%%%%%%%%%%

   %%%%%%%%%%%%%%%%%%%%%%%%%%%%%

\def\lemma#1{\bigskip \par \noindent {\bf #1. Lemma.}}
   %%%%%%%%%%%%%%%%%%%%%%%%%%%%%
\def\cor#1{\bigskip \par \noindent {\bf  #1. Corollary.}}
%%%%%%%%%%%%%%%%%%%%%%%%%%%%%%%%%%%%%%%%%
\def\df#1{\bigskip \par \noindent {\bf  #1. Definition.}}
%%%%%%%%%%%%%%%%%%%%%%%%%%%%%%%%%%%%%%%%%%%%%%
\def\nnt#1{\bigskip \par \noindent {\bf  #1. Notation.}}
%%%%%%%%%%%%%%%%%%%%%%%%%%%%%%%%%%%%%%%%%%%%%%%%%
\def\rem#1{\bigskip \par \noindent {\bf  #1. Remark.}}
%%%%%%%%%%%%%%%%%%%%%%%%%%%%%%%%%%%%%%%%%%%%%%%%%

\vskip 1cm \centerline{\bf  Osculating varieties of Veronesean and
their higher secant varieties.}
\bigskip

\centerline{\it A.Bernardi,
  M.V.Catalisano, A.Gimigliano, M.Id\`a}

\bigskip {\bf Abstract.} We consider the varieties $O_{k,n.d}$ of the k-osculating spaces to the Veronese varieties, the $d-$uple embeddings of $\PP n$; we study the dimension of their higher secant varieties. Via inverse systems (apolarity) and the study of certain spaces of forms we are able, in several cases, to determine whether those secant varieties are defective or not.

\bigskip

\a {\bf  0. Introduction.}

Let us consider the following case of a quite classical problem:
given a generic form $f$ of degree $d$ in $R:=k[x_0,...,x_n]$,
what is the minimum $s$ for which it is possible to write
$f=L^{d-k} _1F_1+...+L^{d-k}_sF_s$, where $L_i\in R_1$ and $F_i\in
R_k$? When $k=0$ this is known as ``Waring problem for forms" (the
original Waring problem is for integers), and it has been solved
via results in \ref {AH}, e.g. see \ref {IK} or \ref {Ge}.

In its generality, this is what was classically called ``to find
canonical forms for a $(n+1)$-ary $d$-ic" (e.g. see \ref {W}).
\par \medskip
We will study this problem here via the study of the dimension of
higher secant varieties to osculating varieties of Veronesean,
since this geometrical problem is equivalent to the one stated
before.

 \a {\bf 1. Preliminaries.}
 \nnt {1.1} \par {\bf i)} In the following we set
$R:=k[x_0,...,x_n]$, where $k=\bar k$ and char$k =0$, hence $R_d$
will denote the forms of degree $d$ on $\PP n$.
\par {\bf ii)} If $X\subseteq \PP N$ is an irreducible projective variety,
an $m$-fat point on $X$ is the $(m-1)^{th}$ infinitesimal
neighborhood of a smooth point $P$ in X, and it will be denoted by
$mP$ (i.e. the scheme $mP$ is defined by the ideal sheaf ${\cal
I}_{P,X}^m \subset {\cal O}_{X}$). 
\b Let dim$X=n$; then, $mP$ is a
0-dimensional scheme of length ${m-1+n\choose n}$. \b If $Z$ is
the union of the $(m-1)^{th}$-infinitesimal neighborhoods in $X$
of $s$ generic points of $X$, we shall say for short that $Z$ is
union of $s$ generic $m$-fat points on $X$.
\par {\bf iii)} If $X\subseteq \PP N$ is a variety and $P$ is a smooth point on
it, the projectivized tangent space to $X$ at $P$ is denoted by
$T_{X,P}$.
\par {\bf iv)} We denote by $<U,V>$ both the linear span in a vector space
or in a projective space of two linear subspaces $U,V$.
\par {\bf v)} If $X$ is a 0-dimensional scheme, we denote by $l(X)$ its
length, while its support is denoted by ${\rm supp} X$.

\a \df {1.2} Let $X\subseteq \PP N$ be a closed irreducible
projective variety;  the $(s-1)^{th}$ {\it higher secant variety}
of $X$ is the closure of the union of all linear spaces spanned by
$s$ points of $X$, and it will be denoted by $X^s$. \b Let $\dim
X=n$; the {\it expected dimension} for $X^s$ is
$$ {\rm expdim} X^s := \hbox{min}\{N, sn+s-1\}$$
where the number $sn+s-1$ corresponds to $\infty ^{sn}$ choices of
$s$ points on $X$, plus $\infty ^{s-1}$ choices of a point on the
$\PP {s-1}$ spanned by the $s$ points. When this number is too
big,  we expect that $X^s = \PP N$. Since it is not always the
case that $X^s$ has the expected dimension, when $ \dim X^s <
\hbox{ min} \{N, sn+s-1\}$, $X^s$ is said to be {\it defective}.

\a A classical result about secant varieties is Terracini's Lemma
(see \ref  {Te}, or, e.g. \ref {A}), which we give here in the
following form:

\a {\bf {1.3}. Terracini's Lemma:} {\it Let $X$  be an irreducible
variety in $\PP N$, and let $P_1,...,P_s$  be s generic points on
$X$. Then, the projectivised tangent space to $X^s$ at a generic
point $Q\in <P_1,...,P_s>$  is the linear span in $\PP N$ of the
tangent spaces $T_{X, P_i}$ to $X$ at $P_i$, $i=1,...,s$, hence}
$$ \dim X^s = \dim <T_{X,P_1},...,T_{X,P_s}>.$$

\cor {1.4} {\it Let $(X,{\cal L})$  be an integral, polarized
scheme.  If ${\cal L}$  embeds $X$ as a closed scheme in $\PP N$,
then
$$ \dim X^s = N - \dim h^0({\cal I}_
{Z,X}\otimes {\cal L}) $$ where Z is union of s generic $2$-fat
points in X.}

\proof By Terracini's Lemma, $ \dim X^s =\dim <T_{X,P_1},...,T_
{X,P_s}> $, with $P_1,...,P_s$ generic points on $X$. Since $X$ is
embedded in $\PP N = \PPP (H^0(X,{\cal L})^*)$, we can view the
elements of $ H^0(X,{\cal L})$ as hyperplanes in $\PP N$; the
hyperplanes which contain a space $T_{X,P_i}$ correspond to
elements in $H^0({\cal I}_ {2P_i,X}\otimes {\cal L})$, since they
intersect $X$ in a subscheme containing the first infinitesimal
neighborhood of $P_i$. Hence the hyperplanes of $\PP N$ containing
the subspace $<T_{X,P_1},...,T_{X,P_s}> $ are the sections of
$H^0({\cal I}_{Z,X}\otimes {\cal L})$, where $Z$ is the scheme
union of the first infinitesimal neighborhoods in $X$ of the
points $P_i$'s. \qquad \qquad \qed

\df {1.5}  Let $X \subset \PP N$ be a variety, and let $P\in X$ be
a smooth point; we define {\it the $k^{th}$ osculating space to X
at P} as the linear space generated by $(k+1)P$, and we denote it
by $O_{k,X,P}$; hence $O_{0,X,P}=\{P\}$, and $O_{1,X,P}=T_{X,P}$,
the projectivised tangent space to $X$ at $P$.
\smallskip \b
Let $X_0\subset X$ be the dense set of the smooth points where
$O_{k,X,P}$ has maximal dimension. The {\it $k^{th}$ osculating
variety to X } is defined as:
$$ O_{k,X}= \overline {\bigcup _{P\in X_0}O_{k,X,P}}. $$

\a \a {\bf 2. Osculating varieties to Veronesean, and their higher
secant varieties.}

\a {\bf  2.1. Notation.} \par {\bf i)} We will consider here
Veronese varieties, i.e. embeddings of $\PP n$  defined by the
linear system of all forms of a given degree $d$: $\nu_d: \PP n
\rightarrow \PP N$, where $N={n+d\choose n}-1$. The $d$-ple
Veronese embedding of $\PP n$, i.e. $Im \nu_d$, will be denoted by
$\, X_{n,d}$.
\par {\bf ii)} In the following we set $O_{k,n,d}:=O_{k,X_{n,d}}$, so that the
$(s-1)^ {th}$ higher secant variety to the $k^{th}$ osculating
variety to the Veronese  variety $X_{n,d}$ will be denoted by
$O_{k,n,d}^s$.

\a \rem {2.2} From now on $\PP N = \Bbb P(R_d)$; a form $M$ will
denote, depending on the situation, a vector in $R_d$ or a point
in $\PP N$.

\b We can view $X_{n,d}$ as given by the map $(\PP n)^{\ast}
\rightarrow \PP N$, where $L \rightarrow L^d$, $L\in R_1$. Hence
$$X_{n,d}=\{ L^d, \quad L\in R_1\}.$$
Let us assume (and from now on this assumption will be implicit)
that $d \geq k$; at the point $P=L^d$ we have (see \ref {Se}, \ref
{CGG} sec.1, \ref {BF} sec.2):

$$O_{k,X_{n,d},P}=\{ L^{d-k}F, \quad F\in R_k\}. \qquad
\qquad (*)$$

\b Notice that $O_{k,X_{n,d},P}$ has maximal dimension $\dim
R_k-1= {{k+n\choose n}-1}$ for all $P\in X_{n,d}$. This can also
be seen in the following way: the fat point $(k+1)P$ on $X_{n,d}$
gives independent conditions to the hyperplanes of $\PP N$, since
it gives independent conditions to the forms of degree $d$ in $\PP
n$. \b Hence, $O_{k,n,d}= \bigcup _ {P\in
X_{n,d}}O_{k,X_{n,d},P}$.

\a As we have already noticed, for $k=0\,$ $(*)$ gives $\,
O_{k,X_{n,d},P}=\{P\}=\{L^d\}\,$, and for $k=1\,$ it becomes $\,
O_{k,X_{n,d},P}=T_{X_ {n,d},P}=\{ L^{d-1}F,\quad  F\in R_1\}$. \a
In general, we have:
$$O_{k,n,d}= \{ L^{d-k}F,  \quad  L\in R_1, \quad F\in R_k\}.$$ \a
Hence,
$$O_{k,n,d}^s=\{L^{d-k}_1F_1+...+L^{d-k}_sF_s, \quad L_i\in R_1,
\quad F_i\in R_k, \quad i=1,...,s\}. $$

\a In the following we also need to know the tangent space
$T_{O_{k,n,d},Q}$ of $O_{k,n,d}$ at the generic point $Q=L^{d-k}F$
(with $L\in R_1, \quad F\in R_k\,$); one has that the affine cone
over $T_{O_{k,n,d},Q}$ is $W=W(L,F)=<L^ {d-k}R_k,L^{d-k-1}FR_1>$
(see \ref {CGG} sec.1, \ref {BF} sec.2)).

\a \lemma  {2.3} The dimension of $O_{k,n,d}$ is always the
expected one, that is, $${\dim} O_{k,n,d} = {\rm min} \{N,\;
n+{k+n\choose n}-1\}$$

\proof   By  2.2, $\dim O_{k,n,d}=\dim W(L,F)-1$, for a generic
choice of $L,F$, so that we can assume that $L$ does not divide
$F$. When $\PPP (W) \neq \PP N$, we have $\dim W = \dim
L^{d-k}R_k+\dim L^{d-k-1}FR_1 - \dim L^{d-k}R_k \cap L^{d-k-1}FR_1
= {k+n\choose n}+(n+1)-1 = {k+n\choose n}+n$, since there is only
the obvious relation between $LR_k$ and $FR_1$, namely $LF-FL=0$.

\a {\bf 2.4.} Consider the classic Waring problem for forms, i.e.
``if we want to write a form of degree $d$ as a sum of powers of
linear forms, how many of them are necessary?" The problem is
completely solved. In fact, $X_{n,d}^s=\{L^d_1+...+L^d_s, \quad
L_i\in R_1\}$ (see previous remark), hence the Waring problem is
equivalent to the problem of computing dim$X_{n,d}^s$. \b By
Corollary 1.4 we have that dim$X_{n,d}^s=N-\dim H^0({\cal I}_
{Z,\PP n}\otimes {\cal O}(d))=H(Z,d)-1$, where $Z$ is a scheme of
$s$ generic 2-fat points in $\PP n$, and $H(Z,d)$ is the Hilbert
function of $Z$ in degree $d$. Since $H(Z,d)$ is completely known
(see \ref{AH}), we are done.

\a More generally, one could ask which is the least $s$ such that
a form of degree $d$ can be written as
$L^{d-k}_1F_1+...+L^{d-k}_sF_s$, with $L_i\in R_1 $ and $F_i\in
R_k$ for $i=1,...,s$; since by Remark 2.2 the variety $O_{k,n,d}
^s$ parameterizes exactly the forms in $R_d$ which can be written
in this way, this is equivalent to answering, for each $k,n,d$, to
the following question:

\a \centerline{\it Find the least $s$, for each $k,n,d$, for which
$O_{k,n,d}^s = \PP N$.}

\a  We are interested in a more complete description of the
stratification of the forms of degree $d$ parameterized by those
varieties, namely in answering the following question:

\a \centerline{\it Describe all $s$ for which $O_{k,n,d}^s\,$ is
defective, i.e. for which $\dim O_{k,n,d}^s < {\rm expdim}
O_{k,n,d}^s.$}

\a Notice that, since $d \geq k$, one has $\dim O_{k,n,d}=N$ if
and only if ${d+n \choose n } \leq n +{k+n \choose n }$, hence for
all such $k,n,d$ and for any $s$ we have $\dim O_{k,n,d}^s = {\rm
expdim} O_ {k,n,d}^s=N$. \b So we have to study this problem when
${d+n \choose n}> n +{k+n \choose n }$, $s\geq 2$; it is easy to
check that whenever $n\geq 2$ this condition is equivalent to
$d\geq k+1$; on the other hand the case $n=1$ (osculating
varieties of rational normal curves) can be easily described (all
the $O_{k,1,d}^s$'s have the expected dimension, see next
section), thus the question becomes:

\a {\bf Question Q(k,n,d)}: {\it For all $k,n,d$ such that $d \geq
k+1$, $n\geq 2$, describe all $s$ for which}
$$\dim O_{k,n,d}^s < \min \{\, N,\, s(n+{k+n\choose n}-1)+s-1\}=\min \{\,
{d+n\choose n}-1,\, s{k+n\choose n}+sn-1\}.$$

\rem {2.5} Terracini's Lemma 1.4 says that  $\dim
O_{k,n,d}^s=N-h^0({\cal I} _ X\otimes {\cal O}_{\PP N}(1)),$ where
$X$ is a generic union of 2-fat points on $O_ {k,n,d}$; we are not
able to handle directly the study of $h^0({\cal I} _ X\otimes
{\cal O}_{\PP N}(1))$, nevertheless, Terracini's Lemma 1.3 says
that the tangent space of $O_{k,n,d}^s$ at a generic point of
$<P_1,...,P_s>, \quad P_i \in O_{k,n,d}\,$, is the span of the
tangent spaces of $O_{k,n,d}$ at $P_i$; if $T_{O_{k,n,d},P_i}=\Bbb
P (W_i) $, then
$$\dim O_{k,n,d}^s=\dim <T_{O_{k,n,d},P_1},...,T_{O_{k,n,d},P_s}>
=\dim <W_1,...,W_s>-1$$

\a We want to prove, via Macaulay's theory of ``inverse systems",
(see \ref {I}, \ref {IK}, \ref {Ge}, \ref {CGG}, \ref {BF}) that,
for a single $W_i$, $\dim W_i=N+1-h^0(\PP n,{\cal I}_ Z(d))$ where
$Z=Z(k,n)$ is a certain 0-dimensional scheme that we will analyze
further, and $\dim \, <W_1,...,W_s>=N+1-h^0(\PP n,{\cal I} _
Y(d))$ where $Y=Y(k,n,s)$ is a generic union in $\PP n$ of $s$
0-dimensional schemes isomorphic to $Z$. Hence,
$$\dim O_{k,n,d}^s=\dim\, <W_1,...,W_s>-1=N - h^0(\PP n,{\cal I}_ Y
(d)).$$

\a So, one strategy in order to answer to the question $Q(k,n,d)$
for a given $(k,n,d)$ is the following:

\par  $1^{st}$ step: try to compute directly $\dim\, <W_1,...,W_s>$; if this is not possible, then
\par \medskip
$2^{nd}$ step: use the theory of inverse systems (classically {\it
apolarity}):

 Compute $W^ {\perp}\subset R_d$, with respect to the perfect pairing $\phi :R_d
\times R_d\rightarrow k$, where: \b - $W$ is a vector subspace of
$R_d$, \b - $\phi (f,g):=\Sigma_{I\in A_{n,d}}f_Ig_I$, where
$A_{n,d}:=\{(i_0,...,i_n) \in \Bbb N^{n+1}, \Sigma_j i_j=d \}$,
with any fixed ordering; this gives a monomial basis
$\{x_0^{i_0}\cdot ...\cdot x_n^{i_n} \}$ for the vector space
$R_d$; if $f \in R_d$, $f=\Sigma_{_{i_0,...,i_n \in
A_{n,d}}}f_{i_0,...,i_n}x_0^{i_0}\cdot ... \cdot x_n^{i_n}$, we
write for short $f=\Sigma f_I{\bf x}^I$, with $I=(i_0,...,i_n)$.

\b Then, consider  $I_d:= W^ {\perp }\subset R_d$. It generates an
ideal $(I_d)\subset R$; in this way we define the scheme
$Z(k,n,d)\subset \PP n$ by setting: $I_{Z(k,n,d)}:= (I_d)^{sat}$.
We will show that these schemes do not depend on $d$.

\par $3^{rd}$ step, compute the postulation for a generic union of $s$ schemes
$Z(k,n,d)$ in $\PP n$.
\smallskip

 \b Recall that $[<W_1,...,W_s>]^{\perp} = W_1^{\perp} \cap
...\cap W_s^{\perp}$.

\lemma {2.6} {\it For all $k,n$ and $d\geq k+2$, we have:
\par
$$(k+1)O \subset Z(k,n,d) \subset (k+2)O,$$
\par
where $Z(k,n,d)$ was defined in {\rm 2.5}, and $O={\rm supp}\
Z(k,n,d) \in \PP n $.}

\a {\bf Proof.} Let $W=<L^{d-k}R_k,\, L^{d-k-1}FR_1> \, \subset
R_d$ be the affine cone over $T_{O_{k,n,d},Q}$ at a generic point
$Q=L^{d-k}F$, with $L\in R_1, \quad F\in R_k\,$. Without loss of
generality we can choose $L=x_0$, so that
$W=x_0^{d-k-1}(x_0R_{k}+FR_1)\,$, hence $x_0^{d-k}R_k\subset W
\subset x_0^{d-k-1}R_{k+1}\,$. So for any $(k,n,d)$,
$$(x_0^{d-k-1}R_{k+1})^{\perp}\subset W^{\perp}\subset (x_0^{d-k}R_k)^
{\perp}.\qquad \qquad (**)$$

 Now, denoting by
$\goth p$ the ideal $(x_1,...,x_n)$, we have:

$$(x_0^{d-t}R_{t})^{\perp}= <\{ x_0^{i_0}\cdot ... \cdot x_n^{i_n}
\, \vert \, \Sigma_j i_j=d, i_0 \leq d-t-1\}>= $$
$$<(\goth p^{d})_d, \, x_0(\goth
p^{d-1})_{d-1},...,x_0^{d-t-1}(\goth p^{t+1})_{t+1}>=(\goth
p^{t+1})_d.$$

\b Now let us view everything in $(**)$ as the degree $d$ part of
a homogeneous ideal; we get:

$$(\goth p^{k+2})_d \subset
(I_{Z(k,n,d)})_{_d }\subset (\goth p^{k+1})_d.$$ Let
$(x_1,...,x_n)$ be local coordinates in $\PP n$ around the point
$O=(1,0,...,0)$; the above inclusions give, in terms of
0-dimensional schemes in $\PP n$:
$$(k+1)O \subset Z(k,n,d) \subset (k+2)O. $$

\lemma {2.7} For any $k,n,d$ with  $d\geq k+2$ the length of
$Z=Z(k,n,d)$ is:
 $$l(Z)=dim W={k+n \choose n}+n.$$

\b {\bf Proof.} We have seen that $Z(k,n,d) \subset (k+2)O$, with
$O$ a point in $\PP n$ (notice that this part of the inclusions in
2.6 works also for $d=k+1$); setting $X:=(k+2)O$, $d \geq k+1$
then gives ${d+n \choose n} \geq l(X)={k+1+n \choose n} \geq
l(Z)$. \b We have ($W \neq R_d$ by assumption) $\dim I_d=dim W^
{\perp}={d+n \choose n}-\dim W$, hence if we prove that $\dim I_d=
{d+n \choose n}-l(Z)$, that is, if $Z$ imposes independent
conditions to the forms of degree $d$, the thesis will follow. \b
One $(k+2)$-fat point always imposes independent conditions to the
forms of degree $d \geq k+1$. Since $Z \subset X= (k+2)O$, then 
$h^1 ({\cal I}_ {Z}(d))=0$ immediately follows.

\bigskip
 Now we have seen that our problem can be translated into a problem of studying certain schemes
 $Z(k,n,d)\subset \PP n$; we want to check that actually these
 schemes are the same for all $d\geq k+2$, say $Z(k,n,d)=Z(k,n)$.
\par
\medskip
\lemma {2.8} {\it For any $k,n$ and  $d\geq k+2$,  we have
$Z(k,n,d)=Z(k,n,k+2)$. Henceforth we will denote  $Z(k,n) =
Z(k,n,d)$, for all $d\geq k+2$.}
\medskip
\b {\bf Proof.} By the previous lemmata we already know that
$Z(k,n,d)$ and $Z(k,n,k+2)$ have the same support and the same
length, hence it is enough to show that $Z(k,n,d)\subset
Z(k,n,k+2)$ (as schemes) in order to conclude.
This will be done if we check that $I(Z(k,n,k+2))_d\subset
I(Z(k,n,d))_d$; in fact, since both ideals are generated in
degrees $\leq d$, this will imply that $I(Z(k,n,k+2))_j\subset
I(Z(k,n,d))_j$, $\forall j\geq d$, hence the inclusion will hold
also between the two saturations, implying $Z(k,n,d)\subset
Z(k,n,k+2)$.
\par
Let $f\in I(Z(k,n,k+2))_d$, then $f=h_1g_1+...+h_rg_r$, where
$h_j\in R_{d-k-2}$ and $g_j\in I(Z(k,n,k+2))_{k+2}$; since
$I(Z(k,n,d))_d$ is the perpendicular to
$W=<L^{d-k}R_k,L^{d-k-1}FR_1>$, it is enough to check that $h_jg_j
\in W^\perp$, $j=1,...,r$. Without loss of generality we can
assume $L=x_0$; hence, since $g_j\in <L^2R_k,LFR_1>^{\perp}$, $g_j
= x_0g'+g''$, with $g',g''\in k[x_1,...,x_n]$ and $g'\in
(FR_1)^{\perp}$. It will be enough to prove
$x_0^{i_0}...x_n^{i_n}g_j=x_0^{i_0+1}...x_n^{i_n}g'+x_0^{i_0}...x_n^{i_n}g''\in
W^\perp$, $\forall i_0,...,i_n$ such that $i_0+...+i_n=d-k-2$. It
is clear that $x_0^{i_0}...x_n^{i_n}g''\in W^\perp$, since
$i_0\leq d-k-2$; on the other hand,  $x_0^{i_0+1}...x_n^{i_n}g'\in
(x_0^{d-k}R_k)^{\perp}$ again by looking at the degree of $x_0$,
while $x_0^{i_0+1}...x_n^{i_n}g'\in (x_0^{d-k-1}FR_1)^{\perp}$
since $g'\in (FR_1)^{\perp}$.

\par
\medskip
\medskip
\rem {2.9} From the lemmata above it follows that in order to
study the dimension of ${\cal O}^s_{k,n,d}$, $\forall d\geq k+2$,
we only need to study the postulation of unions of schemes
$Z(k,n)$. For $d=k+1$, we will work directly on $W$, see
Proposition 3.4.
\par
What we got is a sort of ``generalized Terracini" for osculating
varieties to Veronesean, since the formula $\dim  {\cal
O}^s_{k,n,d}= N-h^0({\cal I}_Y(d))$ reduces to the one in
Corollary 1.4 for $k=0$. Instead of studying 2-fat points on
$X^s_{n,d}$ (see Remark 2.5), we can study the schemes $Y\subset
\PP n$.
\par
\bigskip
\nnt {2.10} Let $Y\subset \PP n$ be a 0-dimensional scheme; we say
that $Y$ is {\it regular} in degree $d$, $d\geq 0$, if the
restriction map $\rho : H^0({\cal O}_{\PP n}(d))\rightarrow
H^0({\cal O}_{Y}(d))$ has maximal rank, i.e. if $h^0({\cal
I}_{Y}(d)).h^1({\cal I}_{Y}(d))=0$. We set $exp\ h^0({\cal
I}_{Y}(d)):= max\ \{0,{d+n\choose n}-l(Y)\}$; hence to say that
$Y$ is regular in degree $d$ amounts to saying that $h^0({\cal
I}_{Y}(d))=exp\ h^0({\cal I}_{Y}(d))$.
\par
Since we always have $h^0({\cal I}_{Y}(d))\geq exp\ h^0({\cal
I}_{Y}(d))$, we write
$$ h^0({\cal I}_{Y}(d))=exp\ h^0({\cal I}_{Y}(d))+\delta,$$
where $\delta=\delta (Y,d)$; hence whenever ${d+n\choose
n}-l(Y)\geq 0$, we have $\delta = h^1({\cal I}_{Y}(d))$, while if
${d+n\choose n}-l(Y)\leq 0$, $\delta = {d+n\choose
n}-l(Y)+h^1({\cal I}_{Y}(d))$; in any case, by setting $exp\
h^1({\cal I}_{Y}(d)):= max\ \{0,l(Y)-{d+n\choose n}\}$, we get:
$h^1({\cal I}_{Y}(d))=exp\ h^1({\cal I}_{Y}(d))+\delta$.
\par
\bigskip
\rem {2.11} {\it For any $k,n,d$ such that $d\geq k+1$, let
$Y=Y(k,n,s)\subset \PP n$ be the 0-dimensional scheme defined in
2.5 for $Z=Z(k,n)$, and $\delta = \delta (Y,d)$. Then
$$\dim O_{k,n,d}^s = {\rm
expdim} O_{k,n,d}^s-\delta.$$ In particular, $\dim O_{k,n,d}^s =
{\rm expdim} O_{k,n,d}^s$ if and only if: \a $h^0({\cal I}_
Y(d))=0, \qquad $ when ${d+n \choose n} \leq s{k+n \choose n }
+sn$; \a $h^0({\cal I}_ Y(d))=N+1-l(Y)={d+n \choose n} - s{k+n
\choose n }-sn$ \  (i.e. $h^1({\cal I}_ Y(d))=0)$, \qquad
 when ${d+n \choose n} \geq s{k+n \choose n }+sn.$}
\a  In fact $h^0({\cal I}_ Y(d))=\ ker \rho $ and  $l(Y)=s{k+n
\choose n }+sn$ (lemma 2.7), ${\rm expdim} O_ {k,n,d}^s=\min \{\,
N={d+n\choose n}-1,\, s(n+{k+n\choose n})-1\}$, and $\dim
O_{k,n,d}^s=N - h^0({\cal I}_ Y(d))=N-exp\ h^0({\cal I}_ Y(d))-
\delta$ (see 2.5).

\a {\bf 3. A few results and a conjecture.} \a Let us consider
first the cases where the question $\bf Q(k,n,d)$ has been
answered.
 \a $\bf Q(k,1,d).$ In this case every $O_ {k,1,d}^s$, with $d\geq k+2$, has
the expected dimension; in fact here $Z(k,1)=(k+2)O$, and the
scheme $Y= \{s$ $(k+2)$-fat points$\}\in \PP 1$ is regular in any
degree $d$. Notice that for $d=k+1$ we trivially have $O_
{k,1,k+1}=\PP N$.
\medskip
\b $\bf Q(1,n,d).$ Here the variety $O_ {1,n,d}$ is the tangential
variety to the Veronese $X_{n,d}$. It is shown in \ref {CGG} that
$Z(1,n)$ is a ``$(2,3)- $scheme" (i.e. the intersection in $\PP n$
of a 3-fat point with a double line); this is easy to see, e.g. by
choosing coordinates so that $L=x_0, \; F=x_1$. \b The postulation
of generic unions of such schemes in $\PP n$, and hence the
defectivity of $O^s_{1,n,d}$, has been studied. Moreover, a
conjecture regarding all defective cases is stated there:
\par
\medskip
\noindent {\bf Conjecture} ( \ref {CGG}). $O^s_{1,n,d}$ is not
defective, except in the following cases:
\par
1) for $d=2$ and $n\geq 2s$;
\par
2) for $d=3$ and $n=s=2,3,4$.
\par
\medskip
In \ref {CGG} the conjecture is proved for $s\leq 5$ (any $d,n$),
for $d=2$ (any $s,n$), for $d\geq 3$ and $n\geq s+1$, for $d\geq
4$ and $s=n$, for $s\geq {1\over 3}{n+2\choose 2}+1$. In \ref B,
the conjecture is proved for $n=2,3$ (any $s,d$).

\a $\bf Q(2,2,d).$ In \ref {BF} it is proved that, for any
$(s,d)\neq (2,4)$, $O_{2,2,d}^s$ has the expected dimension.
\par \bigskip
Now we are going to prove some other cases.
\par
\medskip
The following (quite immediate) lemma describes what can be
deduced about the postulation of the scheme $Y$ from information
on fat points:

\a {\bf 3.1 Lemma.} {\it Let $P_1,...,P_s$ be generic points in
$\PP n$, and set $X:= (k+1)P_1\cup ...\cup (k+1)P_s$, $T:=
(k+2)P_1\cup ...\cup (k+2)P_s$. Now let $Z_i$ be a 0-dimensional
scheme supported on $P_i\,$, $(k+1)P_i\subset Z_i\subset
(k+2)P_i$, with $l(Z_i)=l((k+1)P_i)+n$ for each $i=1,...,n$, , and
set $Y:= Z_1\cup ...\cup Z_s$. Then:
\medskip
\b $Y$ is regular in degree $d$ if one of the following a) or b)
holds:
\smallskip \b a)  ${d+n \choose n} \geq s{k+n+1\choose n}\;$, and $h^1({\cal I}_
{T}(d))=0$;
\smallskip \b b)  ${d+n \choose n} \leq s{k+n\choose n}\;$, and $h^0({\cal I}_
{X}(d))=0$.

\a  $Y$ is not regular in degree $d$, with defectivity $\delta$,
if one of the following c) or d) holds:
\smallskip \b c) $h^1({\cal I}_ {X}(d))>exp\; h^1({\cal I}_ Y(d))$; in this case
$\delta \geq h^1({\cal I}_ {X}(d)))- max \{0, l(Y) - {d+n \choose
n} \}$.
\smallskip \b d)  $h^0({\cal I}_ {T}(d))>exp\; h^0({\cal I}_ Y(d))$; in this
case $\delta \geq h^0({\cal I}_ {T}(d))- max \{0, {d+n \choose
n}-l(Y) \}$.}

\proof The statement follows by considering the cohomology of the
exact sequences:
$$0 \rightarrow {\cal I}_ {T}(d)\rightarrow {\cal I}_
{Y}(d)\rightarrow {\cal I}_ {Y,T}(d)\rightarrow 0$$ and
$$0 \rightarrow {\cal I}_
{Y}(d)\rightarrow {\cal I}_ {X}(d)\rightarrow {\cal I}_
{X,Y}(d)\rightarrow 0$$ where we have: $h^1({\cal I}_
{Y,T}(d))=h^1({\cal I}_ {X,Y}(d))=0$ since those two sheaves are
supported on a 0-dimensional scheme.

\lemma {3.2}  {\it Let $s\geq n+1$ and $d<k+1+j{k+1\over n}$, with
$j=2$ when $s\geq n+2$, and $j=1$ for $s=n+1$. Then $O_{k,n,d}^s$
is not defective and $O_{k,n,d}^s=\PP N$.}

\proof Let $Y\subset \PP n$ be as in 2.5; we have to prove that
$h^0({\cal I}_Y(d))=0$ in our hypotheses.
\par
Let $P_1,...,P_s$ be the support of $Y$; we can always choose a
rational normal curve $C\subset \PP n$ containing $n+2$ of the
$P_i$'s (or just all of them if $s=n+1$). For any hypersurface $F$
given by a section of ${\cal I}_Y(d)$, we have that either
$C\subset F$, or $deg (C\cap F)=nd$; hence, if $nd \leq
(k+1)(n+j)$, where $j=2$ when $s\geq n+2$ and $j=1$ for $s=n+1$,
by Bezout we get $C\subset F$. But this is precisely what our
hypothesis on $d$ says, hence $C\subset F$, but since we can
always find a rational normal curve containing $n+3$ points in
$\PP n$, this would imply that any $P\in \PP n$ is on $F$, i.e.
${\cal I}_Y(d) = 0$.

 \a $\bf Q(k,2,k+2)$. The following corollary describes this
 case completely:

\cor {3.3} {\it Assume $d=k+2\;$ and $n=2$. Then, $ O_{k,n,d}^s$
is not defective for $s\geq3$ and $k\geq 1$, and $ O_{k,n,d}^s$ is
defective for $s=2$ and $k\geq 1$.}

\proof By 3.2, $O_{k,2,k+2}^s$ is not defective for $s\geq 3$ and
$d\geq 3$, i.e. $k\geq 2$; the case $k=1$ is already known by \ref
{B}. \b For $s=2$ and $k\geq 1$, let $Y=Y(k,2) \subset \PP 2$ be
the 0-dimensional scheme defined in {\rm 2.5}; it is easy to check
that $exp\; h^0({\cal I}_ Y(d))=exp\;h^0({\cal I}_ {T}(d))=0$, $T$
denoting the generic union of two $(k+2)$-fat points in $\PP 2$.
Since $T$ is not regular in degree $d=k+2$ for any $k\geq 1$, we
conclude by lemma 3.1 d) that $ O_{k,n,k+2}^s$ is defective with
defectivity $\geq h^0({\cal I}_ {T}(d))=1$ (the only section is
given by the $k+2$-ple line through the two points).
\par \medskip

The following results follow from direct computations on $W$.
\par \medskip
\b $\bf Q(k,n,k+1).$ The following proposition describes this case
completely.
\par \medskip
\prop {3.4} {\it If $s\geq 2$ and $d=k+1$ then
 \par A) if $s\leq n-1$ and the expected dimension is $s{k+n\choose n} +sn$,
  then $O^{s}_{d-1,n,d}$ is defective with defect
$\delta=s^2-s+\sum_{h=2}^s(-1)^h {s\choose h} {k-(h-1)+n \choose
n}$;
  \par {B)} if $s\leq n-1$ and the expected dimension is ${d+n\choose n}$ then
  \par {i)} $O^s_{d-1,n,d}$ is defective with defect $\delta={n-s+d\choose
d}-s(n-s+1)$ if $s < {1\over d}{n-s+d\choose d-1}$;
   \par {ii)} $O^s_{d-1,n,d}=\PP N$ (i.e $O^s_{d-1,n,d}$ is regular) if $s\geq {1\over d}{n-s+d\choose d-1}$;

 \par C) if $s\geq n$ then $O^s_{d-1,n,d}= \PP N$.}

\proof \par{A) }  The variety $O^s_{d-1,n,d}$ wouldn't be
defective if the only relations in $W_1+\cdots +W_s$ would be
those we will be able to find in the proof of Proposition 3.5;
what happens here is that there are too many relations: there are
two kinds of them:
\par {1)} $x_iF_j \in <x_iR_k>\cap <F_jR_1>$ for all $i=0,\ldots
,s-1$ and $j=1,\ldots s$. These relations are exactly $s^2$; then
from these we get a defect of $s^2-s$ (because the number of the
allowed relations in order not to get defectivity is $s$);

\par {2)} $x_ix_jF\in<x_iR_k>\cap<x_jR_k>$ where $i\neq j\in\{ 0,\ldots ,s-1 \}$
 and $F\in R_{k-1}$. The defectivity $\delta$ of $O^s_{d-1,n,d}$
will be $\delta=s^2-s+t$ where $t=\sum_{h=2}^s (-1)^h {s\choose h}
{k-(h-1)+n \choose n}$ is the number of independent forms of type
$x_ix_jF$ with $F\in R_{k-1}$. We can observe that $t$ would be
equal to ${s\choose 2}{k-1+n \choose n}$ if for every $F$
belonging to a base of $R_{k-1}$ the forms $x_ix_jF$ were
independent for all $i\neq j\in \{0,\ldots , s-1\}$; but if $s>2$
this is false: consider for example the following three forms
$F_i=x_iG$, $F_j=x_jG$, $F_l=x_lG$ where $G\in R_{k-2}$ then 
$x_ix_jF_l=x_ix_lF_j=x_jx_lF_i$. Now $t$ would be equal to
${s\choose 2}{k-1+n\choose n}-{s\choose 3}{k-2+n\choose n}$ if for
every $G$ belonging to a base of $R_{k-2}$ the forms of type
$x_ix_jx_lG$ were independent for all $i,j,l\in \{ 0,\ldots ,s-1
\}$ with $i\neq j$, $i\neq l$ and $j\neq l$; but, as before, we
can check that if $s>3$, then $t \leq {s\choose 2}{k-1+n\choose
n}-{s\choose 3}{k-2+n\choose n}+{s\choose 4}{k-3+n \choose n}$.
Proceeding in this way we eventually get that $t=\sum_{h=2}^s
(-1)^h {s\choose h} {k-(h-1)+n \choose n}$.
\par
We can conclude that in this case the defect is
$\delta=s^2-s+\sum_{h=2}^s (-1)^h {s\choose h} {k-(h-1)+n \choose
n}$.

\par {B)} If $s{n+d-1\choose n} +ns \geq {n+d\choose n}$ we
expect that $O^s_{d-1,n,d}=\PP N$. We have that $W_1+\dots + W_s =
<x_{0}R_k,\ldots ,x_{s-1}R_k;F_{1}R_1,\ldots ,F_sR_1>$ in
$K[x_0,\ldots ,x_n]_d$. We can suppose $F_i$'s generic for any
$i=1,\ldots ,s$ in $K[x_s,\ldots ,x_{n}]_d$. Then
$O^s_{d-1,n,d}=\PP N$ if and only if $<F_{1}R_1,\ldots
,F_sR_1>\supset K[x_s,\ldots ,x_{n}]_d :=S_d$; we can actually
just consider the vector space $<F_{1}S_1,\ldots ,F_sS_1>$; since
the $F_i$'s are generic, its dimension will be $min\left\{
{n-s+d\choose d} , s(n-s+1) \right\}$ (e.g. see \ref {MMR}); hence we get that
 \par {i)} if $s(n-s+1)< {n-s+d\choose d}$, then $O^s_{d-1,n,d}$ is
 defective. This happens if and only if $s<{1\over d}{n-s+d\choose
 d-1}$. Then the defect is $\delta={n-s+d\choose d}-s(n-s+1)$.
\par {ii)} if $s(n-s+1)\geq {n-s+d\choose d}$, then $O^s_{d-1,n,d}=\PP
 N$ (for example this is always true for $d\geq n$);

\par {C)} It suffices to prove that $O^s_{d-1,n,d}=\PP N$ for $s=n$.
\par
If $s=n$ and $d=k+1$, the subspace $W_1+\cdots +W_s$ can be
written as $<x_0R_k,F_1R_1, \ldots ,x_{n-1}R_k,F_nR_1>$, which
turns out to be equal to $<x_0R_k, \ldots
,x_{n-1}R_k,x_n^{k+1}>=R_{k+1}$ so $O^n_{d-1,n,d}=\PP N$.
\par
\medskip
For $s\leq n+1$, we have several partial results:
\par \medskip

\prop {3.5} {\it If $s\leq n+1$ and $d\geq 2k+1$ then
$O^s_{k,n,d}$ is regular.}

\proof We have to study the dimension of the vector space
$W_1+\cdots +W_s=<L_1^{d-k}R_k,L_1^{d-k-1}F_1R_1,\ldots ,$
$L_s^{d-k}R_k,L_s^{d-k-1}F_sR_1>$, where $L_1,\ldots ,L_s$ are
generic in $R_1$ and $F_1,\ldots ,F_s$ are generic in $R_k$. Since
$s\leq n+1$, without loss of generality we may suppose
$L_i=x_{i-1}$ for $i=1,\ldots ,s$. Since $d\geq 2k+1$, for $\beta
= d-k \geq 0$, the vector space $W_1+\cdots +W_s$ can be written
as $<x_0^{k+\beta +1}R_k, x_0^{k+\beta}F_1R_1,\ldots
,x_{s-1}^{k+\beta +1}R_k,x_{s-1}^{k+\beta}F_sR_1>$. If we show
that for a particular choice of $F_1,\ldots ,F_s\in R_k$ the
dimension of $W_1+\cdots +W_s=expdim(O^s_{k,n,d})+1$ we can
conclude by semi-continuity that $O^s_{k,n,d}$ has the expected
dimension. Let us consider the case
$F_i=x_ix_{i+1}\widetilde{F}_i$ for $i=1,\ldots ,s-2$,
$F_{s-1}=x_{s-1}x_0\widetilde{F}_{s-1}$ and
$F_s=x_0x_1\widetilde{F}_s$, where the $\widetilde{F}_j$'s are
generic forms in $R_{k-2}$, $j=1,\ldots
 ,n+1$. Let $<x_i^{k+\beta +1}R_k>=: A_i$ and
$<x_i^{k+\beta}F_{i+1}R_1>=:A_i'$, $i=0,\ldots ,s-1$; then we get
$A'_i=<x_i^{k+\beta}x_{i+1}x_{i+2}\widetilde{F}_{i+1}R_1>$,
$i=0,\ldots ,s-3$;
$A'_{s-2}=<x_{s-2}^{k+\beta}x_{s-1}x_0\widetilde{F}_{s-1}R_1>$ and
$A'_{s-1}=<x_{s-1}^{k+\beta}x_0x_1\widetilde{F}_sR_1>$. We can
easily notice that $A_i\cap A_j=\{ 0 \} =A_i'\cap A_j'$ for $i\neq
j$ and $A_i\cap A'_i=<x_i^{k+\beta
+1}x_{i+1}x_{i+2}\widetilde{F}_{i+1}>$, $i= 0,\ldots ,s-3$
(analogously if $i=s-2,s-1$). We can conclude that $dim(W_1+\cdots
+W_s)=s{k+n\choose n}+s(n+1)-s$, which is exactly the expected
dimension.

\prop {3.6} {\it If $s \leq n$ and $k+2\leq d \leq 2k$ then
$O^s_{k,n,d}$ is defective with defect $\delta$ such that:
 \par A) $\delta\geq
{n-s+d\choose d}$ if the expected dimension is ${d+n\choose n}$;
 \par B) $\delta\geq {s\choose 2}{2k-d+n\choose n}$ if the expected dimension is
$s{k+n\choose n} +sn$.}

\proof Let $\beta := d-k$; we can rewrite the vector space
$W_1+\cdots +W_s$ as follows:  $<x_0^\beta R_k,x_0^{\beta
-1}F_1R_1,$ $\ldots ,x_{s-1}^\beta R_k,x_{s-1}^{\beta -1}F_sR_1>$.

\par {A)} We can observe that $(K[x_s,\ldots
,x_n])_d\cap (W_1+\cdots +W_s)=\{ 0 \}$, so if we expect that
$O^s_{k,n,d}=\PP N$ we get a defect $\delta\geq {n-s+d\choose d}$.

\par {B)} Suppose now that $s\left[ {k+n\choose n} +n \right] <
{d+n\choose n}$. If $O^s_{k,n,d}$ would have the expected
dimension we would not be able to find more relations among the
$W_i$'s than $x_i^\beta F_{i+1} \in <x_i^\beta R_k>\cap
<x_i^{\beta -1}F_{i+1}R_1>$, for $i=0,\ldots ,s-1$ (as it happens
in Proposition 3.5). But it's easy to see that $x_i^\beta
x_j^\beta F \in <x_i^\beta R_k>\cap <x_j^{\beta}R_k>$ with $i\neq
j$ and $F\in R_{k-\beta}$. We have exactly ${s\choose 2}$ such
terms for any choice of $F\in R_{k-\beta}$. We can also suppose
that the $F_i\in R_k$ that appear in $W_1+\cdots +W_s$ are
different from $x_j^\beta F$ for any $F\in R_{k-\beta}$ and
$j=0,\ldots ,s-1$ because $F_1,\ldots ,F_s$ are generic forms of
$R_k$. Then we can be sure that the form $x_i^\beta x_j^\beta F$
belonging to $<x_i^\beta R_k>\cap<x_j^\beta R_k>$ ins't one of the
$x_i^\beta F_{i+1}$ that belongs to $<x_i^\beta R_k>\cap
<x_i^{\beta -1}F_{i+1}R_1>$. Now $dim(R_{k-\beta })={k-\beta
+n\choose n}$ then we can find ${s\choose 2}{k-\beta +n\choose n}$
independent forms that give defectivity. Then in the case $s\left[
{k+n\choose n} +n \right] < {d+n\choose n}$ we have
$dim(O^s_{k,n,d})\leq expdim -{s\choose 2}{k-\beta +n\choose n} =
expdim - {s\choose 2}{2k-d+n\choose n}$.

\prop {3.7} {\it If $s=n+1$, $k+2 \leq d \leq 2k$ and
$expdim(O^{n+1}_{k,n,d})= (n+1)\left( {k+n\choose n} +n \right)$
then $O^{n+1}_{k,n,d}$ is defective with defect $\delta\geq {n+1
\choose 2}{2k-d+n \choose n}$.}

\proof The proof of this fact is the same as case B) of the
previous proposition.

\prop {3.8} {\it If $s=n+1$, $n\geq {k+2 \over d-k-2}$, $k+2 < d
\leq 2k$ and $expdim(O^{n+1}_{k,n,d})=N$ then $O^{n+1}_{k,n,d}$ is
defective with defect $\delta\geq{(n+1)(d-k-1)-(d+1)\choose n}$.}

\proof If $k+2< d\leq 2k$, then $2< \beta := d-k \leq k$ and
 we have to study the dimension of $W_1+\cdots + W_{n+1}=<x_0^\beta R_k,x_0^{\beta -1}F_1R_1,
 \ldots , x_n^\beta R_k,x_n^{\beta -1}F_{n+1}R_1>$.  If we expect that
$O^{n+1}_{k,n,d}=\PP N$, it suffices to find a form in $R_d$ which
does not belong to $W_1+\cdots + W_{n+1}$. The forms we are
looking for are:
\par {A)} $x_0^{\beta_0}\cdots x_n^{\beta_n}$ with
$\sum_{i=0}^{n}\beta_i=d$ and $0\leq \beta_i \leq \beta-2$ for all
$i\in\{ 0,\ldots ,n \}$, and
\par {B)} $x_0^{\beta_0}\cdots x_n^{\beta_n}$
with  $\sum_{i=0}^{n}\beta_i=d$, at least one $\beta_i=\beta-1$
and each of the others $\beta_j\leq \beta-2$.
\par\noindent
We will count only how many terms we can find in case A) and then
we will conclude that the defectivity will be grater or equal to
this number.

\par {A)} This case is equivalent to find forms of type $x_0^{d-(\gamma_0+k+2)}\cdots
x_n^{d-(\gamma_n+k+2)}$ with
$\sum_{i=0}^{n}\gamma_i=nd-(n+1)(k+2)$ and $\gamma_i\geq 0$ for
all $i=0,\ldots ,n$. Then these forms are exactly
${n+(n+1)(d-k-2)-d\choose n}={(n+1)(d-k-1)-(d+1)\choose n}$. This
will be possible only if $(n+1)(d-k-2)-d\geq 0$ and so if $n\geq
{k+2 \over d-k-2}$.
\par
\bigskip
All the results on defectivity lead us to formulate the following:
\par \medskip
\noindent {\bf 3.9 Conjecture.} {\it $O^{s}_{k,n,d}$ is defective
only if $Y$ is as in case c) or d) of Lemma 3.1.}
\par
\medskip
The conjecture amounts to say that the defectivity of $Y$ can only
occur if defectivity of the fat points schemes $X$ or $T$ imposes
it.
\par \medskip
In a forthcoming paper we intend to explore more in depth the connections between the postulation of fat points and our schemes $Y$.
\par
\bigskip
{\bf Acknowledgments}  All authors supported by MIUR. The last two
authors supported by the University of Bologna, funds for selected
research topics.

\bigskip\bigskip\centerline {{\bf REFERENCES}}

\medskip\noindent[{\bf AH} ]: J. Alexander, A. Hirschowitz. {\it Polynomial
interpolation in several variables.} J. of Alg. Geom. {\bf 4}
(1995), 201-222.

\medskip\noindent [{\bf B}]: E. Ballico, {\it On the
secant varieties to the tangent developable of a Veronese
variety}, Preprint.

\medskip\noindent [{\bf BF}]: E. Ballico, C.Fontanari,{\it On the
secant varieties to the osculating variety of a Veronese surface},
Central Europ. J. of Math. {\bf 1} (2003), 315-326.

\medskip\noindent [{\bf CGG}]: M.V.Catalisano, A.V.Geramita, A.Gimigliano.
{\it On the Secant Varieties to the Tangential Varieties of a
Veronesean.}  Proc. A.M.S.  {\bf 130} (2001), 975-985.

\medskip\noindent [{\bf Ge}]: A.V.Geramita. {\it Inverse Systems of Fat
Points}, Queen's Papers in Pure and Applied Math.  {\bf 102},
{\it The Curves Seminar at Queens', vol. X} (1998).

\medskip\noindent [{\bf IK}]: A.Iarrobino, V.Kanev. {\it Power Sums,
Gorenstein algebras, and determinantal loci.} Lecture Notes in
Math. {\bf 1721}, Springer, Berlin, (1999).

\medskip\noindent [{\bf I}]: A.Iarrobino. {\it Inverse systems of a symbolic
algebra III: Thin algebras and fat points.} Compos. Math. {\bf
108} (1997), 319-356.

\medskip\noindent [{\bf MMR}]: J.Migliore, R.Mir\`o -Roig. {\it Ideals of generic forms and the ubiquity of the weak Lefschetz property.} J. Pure Appl. Algebra  {\bf 182}  (2003), 79-107.

\medskip\noindent [{\bf Se}]: B. Segre, {\it Un'estensione delle variet\`a di
Veronese ed un principio di dualit\`a per le forme algebriche I
and II}. Rend. Acc. Naz. Lincei (8) {\bf 1} (1946), 313-318 and
559-563.

\medskip\noindent [{\bf Te}]: A.Terracini. {\it Sulle} $V_k$ {\it per cui la
variet\`a degli} $S_h$ $(h+1)${\it -seganti ha dimensione minore
dell'ordinario.} Rend. Circ. Mat. Palermo {\bf 31} (1911),
392-396.
\medskip
\noindent [{\bf W}]: K.Wakeford. {\it On canonical forms.} Proc.
London Math. Soc. (2) {\bf 18} (1919/20). 403-410.
\par
\medskip

\bigskip \bigskip

\b {\it A.Bernardi, Dip. Matematica, Univ. di Milano, Italy,
email: bernardi@mat.unimi.it}
 \medskip
\b {\it M.V.Catalisano, DIPEM, Univ. di Genova, Italy,
e-mail: catalisano@dimet.unige.it}
\medskip
\b {\it A.Gimigliano, Dip. di Matematica and C.I.R.A.M., Univ. di
Bologna, Italy, e-mail: gimiglia@dm.unibo.it}
\medskip
\b {\it M.Id\`a, Dip. di Matematica, Univ. di Bologna, Italy,
e-mail: ida@dm.unibo.it}

\end